\documentclass[preprint,authoryear,12pt]{elsarticle}

\usepackage{amssymb}

\usepackage{amsthm}
\usepackage[sumlimits]{amsmath}
\usepackage{amscd}
\usepackage{amsfonts}
\usepackage{color,rotating}
\usepackage[margin=2cm]{geometry}
\usepackage{graphicx}
\usepackage{psfrag}
\psfrag{test}{$P1^3$}

\DeclareMathOperator{\diff}{d}
\DeclareMathOperator{\divv}{div}

\DeclareMathOperator{\Ro}{Ro}
\DeclareMathOperator{\Fr}{Fr}

\newtheorem{theorem}{Theorem}

\newtheorem{definition}[theorem]{Definition}

\def\MM#1{\boldsymbol{#1}}
\newcommand{\pp}[2]{\frac{\partial #1}{\partial #2}} 

\newcommand{\dd}[2]{\frac{\diff#1}{\diff #2}}

\usepackage{xspace}

\newcommand{\pdg}{${P}1_{{DG}}$\xspace}

\newcommand{\verts}{\mathrm{vert}}
\newcommand{\face}{\mathrm{face}}
\newcommand{\edge}{\mathrm{edge}}

\journal{Journal of Computational Physics}

\begin{document}

\begin{frontmatter}



\title{Mixed finite elements for numerical weather prediction}


\author[aero]{C.~J. Cotter}
\author[aero]{J. Shipton}
\address[aero]{Department of Aeronautics, Imperial College London,
 South Kensington Campus, London SW7 2AZ}

\begin{abstract}
  We show how two-dimensional mixed finite element methods that
  satisfy the conditions of finite element exterior calculus can be
  used for the horizontal discretisation of dynamical cores for
  numerical weather prediction on pseudo-uniform grids. This family of
  mixed finite element methods can be thought of in the numerical
  weather prediction context as a generalisation of the popular
  polygonal C-grid finite difference methods. There are a few major
  advantages: the mixed finite element methods do not require an
  orthogonal grid, and they allow a degree of flexibility that can be
  exploited to ensure an appropriate ratio between the velocity and
  pressure degrees of freedom so as to avoid spurious mode branches in
  the numerical dispersion relation. These methods preserve several
  properties of the C-grid method when applied to linear barotropic
  wave propagation, namely: a) energy conservation, b) mass
  conservation, c) no spurious pressure modes, and d) steady
  geostrophic modes on the $f$-plane. We explain how these properties
  are preserved, and describe two examples that can be used on
  pseudo-uniform grids: the recently-developed modified RT0-Q0
  element pair on quadrilaterals and the BDFM1-\pdg element pair on
  triangles. All of these mixed finite element methods have an exact
  2:1 ratio of velocity degrees of freedom to pressure degrees of
  freedom. Finally we illustrate the properties with some numerical
  examples.
\end{abstract}


\end{frontmatter}

\section{Introduction}

There are a number of groups that have been developing dynamical cores
for numerical weather prediction (NWP) and climate modelling, based on
triangular meshes on the sphere or on the dual meshes composed of
hexagons together with twelve pentagons
\citep{Ri+2000,Ma+2002,Sa+2008}. These grids are referred to as
pseudo-uniform grids since they have edge lengths $h$ that satisfy
$c_0\bar{h}<h<c_1\bar{h}$, as $\bar{h}\to 0$, where $\bar{h}$ is the
average edge length, for some positive constants $c_0$, $c_1$. The
principal reason for adopting these grids is that they provide a
direct addressing data structure whilst avoiding the polar singularity
of the latitude-longitude grid, which introduces a bottleneck to
scaling on massively parallel architectures due to the convergence of
meridians. One approach to developing numerical discretisations on
triangular or hexagonal grids is to adapt the staggered Arakawa C-grid
finite difference method on quadrilaterals \citep{ArLa77} (used in
several currently operational NWP models, such as the UK Met Office
Unified Model \citep{Da+2005}) since this type of staggering prevents
pressure modes (non-constant functions on the pressure grid that have
zero numerical gradient). By defining discrete curl and divergence
operators which satisfy div curl$=0$, it is possible
to construct C-grid discretisations for horizontal wave propagation
which have stationary geostrophic modes on the $f$-plane
\citep{ThRiSkKl2009}, a necessary condition for accurate
representation of geostrophic adjustment processes. These operators
can be used to construct energy and enstrophy C-grid discretisations
for the nonlinear rotating shallow-water equations using the vector
invariant form \citep{RiThKlSk2010}. The drawback with using the
C-grid finite difference method on triangles or hexagons instead of
quadrilaterals is that the ratio of velocity and pressure degrees of
freedom (DOF) is altered. The quadrilateral C-grid has one pressure
DOF stored at the centre of each grid cell, and two velocity DOF per
grid cell (normal velocity is stored at each of the four edges, which
are each shared with the neighbouring cell on the other side of the
face)\footnote{Here, and in the rest of the paper, we consider compact
  domains without boundary such as the sphere and rectangles with
  double periodic boundary conditions.}. This is considered the ideal
ratio, since the velocity then has an equal number of rotational and
divergent DOF which are coupled together in the correct way so that
there are two inertia-gravity modes (the inward and outward
propagating modes) for each Rossby mode. On the other hand, the
triangular C-grid has only $3/2$ velocity DOF per grid cell, and the
hexagonal C-grid has $3$ velocity DOF per grid cell. This means that
the triangular C-grid has four inertia-gravity modes per Rossby mode;
the extra spurious inertia-gravity branch has a frequency range that
decreases with Rossby deformation radius, leading to ``checkerboard
patterns'' in the divergence when the deformation radius is small (as
it can be in the ocean, or when there are many vertical layers). The
hexagonal C-grid has an equal number of inertia-gravity and Rossby
modes; the extra spurious Rossby mode has very low frequencies and
propagates Eastwards on the $\beta$-plane \citep{Th08}. The effects of
these spurious Rossby modes has not been reported in practice but
there are concerns amongst the operational NWP community that if
spurious modes are supported by the grid, then they might be
initialised during the data assimilation process or by physics
parameterisations \citep{St_pc2010}. It may also be the case that the
spurious modes lead to spurious spread/lack of spread in ensemble
forecasts. Careful numerical experiments are required to investigate
this concern.

The finite element method provides the opportunity to alter the number
of degrees of freedom per triangular element to ameliorate this
problem. A number of finite element pairs on triangles have been
proposed for geophysical fluid dynamics, mostly in the ocean modelling
community \citep{WaCa1998,RoStLi1998,Ro+2005,CoHaPa2009,CoLaReLe2010,
  RoRoPo2007}. In \citep{RoRo2008}, the lowest order
Brezzi-Douglas-Marini element pair \citep{BrDoMa1985}, known as BDM1,
was investigated in the context of the discrete shallow-water
equations. The velocity space is piecewise linear with continuous
normal components, and the pressure space is piecewise constant. The
natural data structure for the velocity space stores two normal
velocity components on each edge, and hence there are 3 velocity DOF
per triangular element and 1 pressure DOF. There are too many velocity
DOF and hence there will be too many Rossby modes per inertia-gravity
mode, just as for the hexagonal C-grid. 

The key result of this paper is in showing that discretisations of the
linear rotating shallow water equations on the $f$-plane constructed
using these spaces on arbitrary meshes satisfy a crucial property,
namely that geostrophic modes are exactly steady. This is achieved by
making use of the discrete Helmholtz decomposition, within the
framework of discrete exterior calculus \citep{ArFaWi2006}. As
described in \citep{Ar2002}, existence of such a decomposition
requires that the following diagram commutes:
 \begin{equation}
 \begin{CD}
 \label{big commutes}
 H^1(\Omega) @>\nabla^\perp>> H(\mathrm{div},\Omega) @>\nabla\cdot>> L_2(\Omega) \\ 
 @VV{\Pi^E}V @VV{\Pi^S}V @VV{\Pi^V}V \\
 E@>\nabla^\perp>> S @>\nabla\cdot>> V
 \end{CD}
 \end{equation}
 where $\Pi^E$, $\Pi^S$ and $\Pi^V$ are suitably chosen projection
 operators.  The same Helmholtz decomposition can then be used to
 study the discrete dispersion relations for the numerical
 discretisation. Within this framework, we then conclude that an
 optimal choice is to have $\dim(S)=2\dim(V)$ which, at least in the
 periodic plane, satisfies necessary conditions for absence of both
 spurious inertia-gravity and spurious Rossby waves.

The rest of this paper is organised as follows. The general framework
of mixed finite element methods applied to the linear rotating
shallow-water equations is described in Section \ref{framework}, and
the four properties of energy conservation, local mass conservation,
absence of spurious pressure modes and steady geostrophic modes are
discussed. In Section \ref{examples}, two examples are then introduced
that fit into this framework, and numerical results are presented in
section \ref{numerics}. Finally, we give a summary and outlook in
Section \ref{summary}.

\section{Mixed finite elements for geophysical fluid dynamics}
\label{framework}
In this section we describe how mixed finite elements can be used to
build flexible discretisations on pseudo-uniform grids. We concentrate
on the rotating shallow-water equations which are regarded in the
numerical weather prediction community as being a simplified model
that contains many of the issues arising in the horizontal
discretisation for dynamical cores. Since in this paper we are
concerned with wave propagation properties, we restrict attention to
the linearised equations on the $f$-plane, $\beta$-plane or the
sphere. First, we introduce the mixed finite element formulation
applied to the linear rotating shallow-water equations, then we
discuss various properties of the formulation that are a requirement
for numerical weather prediction applications, namely global energy
and local mass conservation, absence of spurious pressure modes and
steady geostrophic states. These properties all rely on exact sequence
properties, \emph{i.e.} div-curl relations, as described in
\citep{ArFaWi2006}.

\subsection{Spatial discretisation of the linear rotating shallow-water equations}
In this paper we consider the discretisation of the linear rotating
shallow-water equations on a two dimensional surface $\Omega$ that is
embedded in three dimensions (which we restrict to be compact with no
boundaries, \emph{e.g.} the sphere or double periodic $x-y$ plane):
\begin{equation}
\label{swe}
\MM{u}_t + f\MM{u}^\perp + c^2\nabla\eta = 0, \quad
\eta_t + \nabla\cdot\MM{u} = 0, \quad \MM{u}\cdot\MM{n}=0 \quad
\mathrm{on}\quad \partial\Omega,
\end{equation}
where $\MM{u}=(u,v)$ is the horizontal velocity,
$\MM{u}^\perp=\MM{k}\times\MM{u}$, $f$ is the Coriolis parameter,
$c^2=gH$, $g$ is the gravitational acceleration, $H$ is the mean layer
thickness, $h=H(1+\eta)$ is the layer thickness, $\MM{k}$ is the
normal to the surface $\Omega$, and $\nabla$ and $\nabla\cdot$ are
appropriate invariant gradient and divergence operators defined on the
surface. We form the finite element approximation by multiplying by
time-independent test functions $\MM{w}$ and $\phi$, integrating over
the domain, integrating the pressure gradient term $c^2\nabla\eta$ by
parts in the momentum equation, and finally restricting the velocity
trial and test functions $\MM{u}$ and $\MM{w}$ to a finite element
subspace $S\subset H(\divv)$ (where $H(\divv)$ is the space of square
integrable velocity fields whose divergence is also square
integrable), and the elevation trial and test functions $\eta$ and
$\alpha$ to the finite element subspace $V\subset L^2$ (where $L^2$ is
the space of square integrable functions):
\begin{eqnarray}
\label{fe u eqn}
\dd{}{t}\int_{\Omega}\MM{w}^h\cdot\MM{u}^h\diff{V} +
\int_{\Omega}f\MM{w}^h\cdot\left(\MM{u}^h\right)^\perp\diff{V} -
c^2\int_{\Omega}\nabla\cdot\MM{w}^h\eta^h\diff{V} &=& 0, \quad
\forall\MM{w}^h\in S, \\
\dd{}{t}\int_{\Omega}\alpha^h\eta^h\diff{V} +
\int_{\Omega}\alpha^h\nabla\cdot\MM{u}^h\diff{V}&=&0,
\quad \forall \alpha^h \in V.
\label{fe eta eqn}
\end{eqnarray}
After discretisation in time, these equations are solved in practise
by introducing basis expansions for $\MM{w}^h$, $\MM{u}^h$, $\eta^h$,
and $\alpha^h$ and solving the resulting matrix-vector systems for the
basis coefficients.

In this framework we restrict the choice of finite element spaces $S$
and $V$ so that
\[
\MM{u}^h \in S \quad \implies \nabla\cdot\MM{u}^h \in V.
\]
The divergence should map from S onto V, so that for all functions
$\phi^h\in V$ there exists a velocity field $\MM{u}^h\in S$ with
$\nabla\cdot\MM{u}^h=\phi^h$. Such spaces are known as
``div-conforming''. Furthermore we require that there exists a
``streamfunction'' space $E\subset H^1$ such that
\[
\psi^h \in E \quad \implies \MM{k}\times\nabla\psi^h \in S, 
\]
where the $\MM{k}\times\nabla$ operator (the curl, which we shall
write as $\nabla^\perp$) maps onto the kernel of $\nabla\cdot$ in S. A
consequence of these properties is that functions in $E$ are
continuous, vector fields in $S$ only have continuous normal
components and functions in $V$ are discontinuous.

\subsection{Energy conservation}
Global energy conservation for the linearised equations is a
requirement of numerical weather prediction models for various
reasons, in particular because it helps to prevent numerical sources
of unbalanced fast waves. It is also a precursor to a
energy-conserving discretisation of the nonlinear equations using the
vector-invariant formulation. For the mixed finite element method, global
energy conservation is an immediate consequence of the Galerkin finite
element formulation. The conserved energy of equations \eqref{swe} is
\[
H = \frac{1}{2}\int_{\Omega} |\MM{u}|^2 + c^2 \eta^2 \diff{V}.
\]
Substituting the solutions $\MM{u}^h$ and $\eta^h$ to
equations (\ref{fe u eqn}-\ref{fe eta eqn}) and taking the time
derivative gives
\[
\dd{}{t}H = \int_{\Omega}\MM{u}^h\cdot\dot{\MM{u}}^h + c^2\eta^h\dot{\eta}^h\diff{V}.
\]
Choosing $\MM{w}^h=\MM{u}^h$ and $\alpha^h=\eta^h$ 
in equations (\ref{fe u eqn}-\ref{fe eta eqn}) then gives
\begin{eqnarray*}
\dd{}{t}H &=&
\int_{\Omega}\MM{u}^h \cdot\dot{\MM{u}}^h 
+ c^2\eta^h\dot{\eta}^h\diff{V} \\
& = & \int_{\Omega}-f\underbrace{\MM{u}^h
\cdot\left({\MM{u}}^h\right)^\perp}_{=0} 
+ \underbrace{c^2\nabla\cdot\MM{u}^h\eta^h
- c^2\eta^h\nabla\cdot\MM{u}^h}_{=0}\diff{V} = 0.
\end{eqnarray*}

\subsection{Local mass conservation}
Local mass conservation is a requirement for numerical weather
prediction models since it prevents spurious sources and sinks of
mass. For the nonlinear density equation, this can be achieved using a
finite volume or discontinuous Galerkin method. For mixed finite
element methods of the type used in this paper applied to the linear
equations, consistency and discontinuity of functions in $V$ requires
that element indicator functions (\emph{i.e.} functions that are equal
to 1 in one element and 0 in the others) are contained in
$V$. Selecting the element indicator function for element $e$ as the
test function $\alpha^h$ in equation \eqref{fe eta eqn} gives
\[
\dd{}{t}\int_e \eta^h\diff{V} + \int_{\partial e}
\MM{u}^h\cdot\MM{n}\diff{S} = 0,
\]
where $\partial e$ is the boundary of element $e$. Since
$\MM{u}^h$ has continuous normal components on element
boundaries, this means that the flux of $\eta^h$ is continuous
and hence $\eta^h$ is locally conserved.

\subsection{Absence of spurious pressure modes and stability of discrete Poisson equation}
\label{pressure modes}
The principle reason for using the staggered C-grid for numerical
weather prediction is that the collocated A-grid, in which pressure
and both components of velocity are stored at the same grid locations,
suffers from a checkerboard pressure mode which has vanishing
numerical gradient when the centred difference approximation is used,
despite being oscillatory in space. This pressure mode rapidly
pollutes the numerical solution in the presence of nonlinearity,
boundary conditions and forcing, and can be easily excited by physics
subgrid parameterisations or initialisation \emph{via} data
assimilation from noisy data.

In the context of mixed finite element methods applied to the equation
set \eqref{swe}, spurious pressure modes relate to the discretised gradient
$D\phi^h\in S$ of a function $\phi^h\in V$ defined by
\[
\int_{\Omega}\MM{w}^h\cdot D\phi^h\diff{V} =
-\int_{\Omega}\nabla\cdot\MM{w}^h\phi^h\diff{V}, \quad \forall \MM{w}^h\in S.
\]
On uniform grids, spurious pressure modes are functions $\phi^h$ from
the pressure space $V$ which have zero discretised gradient
$D\phi^h$ even though $\nabla\phi^h$ is non-zero.  On unstructured
grids or grids with varying edge lengths, spurious pressure modes are
functions which have discretised gradient becoming arbitrarily small
as the maximum edge length $h_0$ tends to zero, despite their actual
gradient staying bounded away from zero. Such functions would prevent
the numerical solution of equations \eqref{swe} converging at the
optimal rate predicted by approximation theory. We make the following 
definition of a spurious pressure mode.
\begin{definition}[Spurious pressure modes]
  A mixed finite element space $(S,V)$ is said to be free of spurious
  pressure modes if there exists $\gamma_2>0$ independent of $h_0$
  such that for all $\phi^h\in V$, there exists nonzero $\MM{v}^h\in
  S$ satisfying
\begin{equation}
\label{inf sup}
\int_{\Omega}\phi^h\nabla\cdot\MM{v}^h\diff{V}\geq \gamma_2\|\phi^h\|_{L_2}
\|\MM{v}^h\|_{H(\mathrm{div})}.
\end{equation}
\end{definition}
Condition \eqref{inf sup} is one of two sufficient conditions for
numerical stability of the mixed finite element discretisation of the
Poisson equation $-\nabla^2\phi=f$ given by
\begin{eqnarray*}
\int_{\Omega}\MM{w}^h\cdot \MM{v}^h\diff{V} &=&
-\int_{\Omega}\nabla\cdot\MM{w}^h\phi^h\diff{V}, \quad \forall \MM{w}^h\in S, \\
\int_\Omega \alpha^h \nabla\cdot\MM{v}^h\diff{V} & = & 
\int_\Omega \alpha^h f^h, \quad \forall \alpha^h \in V.
\end{eqnarray*}
This discretisation is stable (\emph{i.e.} small changes in the
right-hand side lead to small changes in the solution field in the
limit as $h_0$ tends to zero) if Condition \eqref{inf sup} holds,
together with the condition that there exists $\gamma_1>0$ independent
of $h_0$ such that
\begin{equation}
\label{S1}
\int_{\Omega}\MM{v}^h\cdot\MM{v}^h\diff{x} \geq \gamma_1 
\|\MM{v}^h\|^2_{H(\mathrm{div})},
\end{equation}
for all $\MM{v}^h\in S$ such that
$\int\nabla\cdot\MM{v}^h\phi^h\diff{V}=0$ for all $\phi^h\in V$.
As reviewed in \citet{Ar2002}, Condition \eqref{inf sup} is satisfied if it is
possible to define a bounded projection $\Pi^S:H(\mathrm{div})\to S$
such that the following diagram commutes:
\begin{equation}
\label{commutes}
\begin{CD}
 H(\mathrm{div},\Omega) @>\nabla\cdot>> L_2(\Omega) \\ 
@VV{\Pi^S}V @VV{\Pi^V}V \\
 S @>\nabla\cdot>> V
\end{CD}
\end{equation}
where $\Pi^V$ is the usual $L_2$ projection operator. This means that
taking any square integrable velocity field $\MM{u}$ with square
integrable divergence, evaluating the divergence and projecting into
$V$ produces the same result as projecting $\MM{u}$ into $S$ using
$\Pi^S$ and evaluating the divergence. The projection $\Pi^S$ is
constructed by applying an $L^2$ projection of normal components on
element edges, ensuring that $\MM{u}$ is $L^2$-orthogonal to gradients
of functions from $V$ in each element, and ensuring the remaining
degrees of freedom in $\MM{u}$ are $L^2$-orthogonal to divergence-free
functions in each element. We shall explain how this is done for the
two examples described in Section \ref{examples}.  To check that the
diagram \eqref{commutes} commutes, it is sufficient to show that
\[
\int_K \alpha^h (\nabla\cdot\MM{u}-\nabla\cdot\Pi^S\MM{u})\diff{V}=0,
\quad \forall \alpha^h\in V, \MM{u}\in H(\mathrm{div},K),
\]
for each element $K$, since this defines the $L_2$ projection $\Pi^V$
into the discontinuous space $V$. This is easily checked using integration by parts:
\begin{eqnarray*}
\int_K\alpha^h\nabla\cdot\MM{u}\diff{V} & = & -
\int_K\nabla\alpha^h\cdot\MM{u}\diff{V} + \int_{\partial
K}\alpha^h\MM{u}\cdot\MM{n}\diff{S}, \\
& = & -
\int_K\nabla\alpha^h\cdot\Pi^S\MM{u}\diff{V} + \int_{\partial
K}\alpha^h\Pi^S\MM{u}\cdot\MM{n}\diff{S}
=  \int_K\alpha^h\nabla\cdot\Pi^S\MM{u}\diff{V},
\end{eqnarray*}
as required.

As also reviewed in \citet{Ar2002}, Condition \eqref{S1} is satisfied
if vector fields $\MM{v}\in S$ with divergence orthogonal to $V$ are
in fact divergence-free. This is satisfied by the types of mixed
finite element methods considered in this paper since the divergence
maps from S into V, and so the projection of $\nabla\cdot\MM{v}^h$
into $V$ is simply the inclusion. Hence, if the divergence is
orthogonal to $V$, the divergence must be zero, and so \eqref{S1} is
satisfied.

\subsection{Discrete Helmholtz decomposition}
Proof of the condition that geostrophic modes are steady requires the
construction of a discrete Helmholtz decomposition.  Since Condition
\eqref{inf sup} holds, the discrete gradient operator $D:V\to S$, has no
non-trivial kernel. For any $\psi^h\in E$, the curl $\nabla^\perp$ of
$\psi^h$ satisfies
\[
\int_\Omega \nabla^\perp\psi^h\cdot D\phi^h\diff{V}=
-\int_\Omega \underbrace{\nabla\cdot\nabla^\perp\psi^h}_{=0}\phi^h\diff{V}=0,
\]
for any $\phi^h\in V$, and hence the curl from $E$ to $S$ and the
discrete divergence from $V$ to $S$ map onto orthogonal subspaces of
$S$.  This means that there is a one-to-one mapping between elements
of $S$ and $E\times V$, defining a discrete Helmholtz decomposition
\begin{equation}
\label{helmholtz}
\MM{u}^h = \nabla^\perp\psi^h + D\phi^h + \MM{h}^h, \quad \MM{u}\in S,\,
\psi^h\in E,\, \phi^h\in V, \, \MM{h}^h\in H, 
\end{equation}
where $H\subset S$ is the space of discrete harmonic velocity fields
\[
H^h = \left\{\MM{u}^h\in S: \nabla\cdot\MM{u}^h=0, \quad
\int_{\Omega} \MM{u}^h\cdot\nabla^\perp\psi^h\diff{V}=0, \,
\forall \psi^h\in E\right\}.
\]
The dimension of $H^h$ is the same as the dimension of the space
$H$ of harmonic velocity fields
\[
H = \left\{\MM{u} \in H(\divv):\nabla\cdot\MM{u}=0, \qquad \int_{\Omega} \MM{u}^h\cdot\nabla^\perp\psi\diff{V}=0, \,\forall \psi\in H^1\right\},
\]
\emph{i.e.}, velocity fields with vanishing divergence and (weak) curl
(In the periodic plane, these harmonic velocity fields are the
constant velocity fields, but there are no harmonic velocity fields on
the sphere); however $H^h \ne H$ in the general case
\citep{ArFaWi2006}. The kernel of $\nabla^\perp$ in $E$ is the
subspace of constant functions, and stability results (as described in
Section \ref{pressure modes}) imply that the kernel of $D$ in $V$ is
the subspace of constant functions, and hence we can use Equation
\eqref{helmholtz} to obtain a DOF count for $S$.
\[
\dim(S)=(\dim(E)-1)+(\dim(V)-1)+\dim(H),
\]
and hence
\[
\dim(E) = \dim(S)-\dim(V)+2-\dim(H).
\]
For our DOF requirement $\dim(S)=2\dim(V)$, we obtain
\[
\dim(E) = \dim(V) + 2 - \dim(H),
\]
which becomes $\dim(E)=\dim(V)$ for the periodic plane and
$\dim(E)=\dim(V)+2$ for the sphere. If $\dim(S)>2\dim(V)$, then
$\dim(E)>\dim(V)+(2-\dim(H))$ and \emph{vice versa}. This will become
important when we examine wave propagation in
Section \ref{dispersion relations}.

\subsection{Vorticity and divergence}

The discrete vorticity associated with the velocity $\MM{u}^h\in S$ is
defined as $\xi^h\in E$ such that
\begin{equation}
\label{vorticity}
\int_{\Omega} \gamma^h\xi^h\diff{V} =
-\int_{\Omega}\nabla^\perp\gamma^h\cdot\MM{u}^h\diff{V},
\quad \forall\gamma^h\in E.
\end{equation}
It is possible to obtain $\MM{u}\in S$ from the discrete vorticity
$\xi\in E$ and the divergence $\delta^h=\nabla\cdot\MM{u}^h\in V$ by
solving two elliptic problems for the streamfunction $\psi^h$ and
velocity potential $\phi^h$. To obtain the streamfunction $\psi^h\in E$,
we use the Helmholtz decomposition and rewrite equation
\eqref{vorticity} as
\[
\int_{\Omega} \gamma^h\xi^h\diff{V} =
-\int_{\Omega}\nabla\gamma^h\cdot\nabla\psi^h\diff{V}, \quad
\forall\gamma^h\in \left\{\gamma:\gamma\in E, \quad
\int_\Omega\gamma\diff{V}=0\right\}, \qquad \int_\Omega \psi^h\diff{V}=0,
\]
which is the usual finite element discretisation of the Poisson
equation for $\psi^h$. To obtain the vector potential $\phi^h$
requires the solution of the coupled system
\begin{eqnarray*}
\int_{\Omega}\alpha^h\nabla\cdot D\phi^h\diff{V} & = & \int_{\Omega}
\alpha^h\delta^h\diff{V}, \quad \forall\alpha^h\in \left\{
\alpha:\alpha\in V,\quad \int_\Omega\alpha\diff{V}=0\right\}, 
\\
\int_{\Omega}\MM{w}^h\cdot D\phi^h\diff{V} 
& = & -\int_{\Omega}\nabla\cdot\MM{w}^h\phi^h\diff{V}, \quad
\forall\MM{w}^h\in S, \qquad \int_\Omega\phi^h\diff{V}=0.
\end{eqnarray*}
This is the mixed finite element approximation to the Poisson equation
already discussed in Section \ref{pressure modes}; if the Conditions
\eqref{inf sup} and \eqref{S1} are satisfied, the coupled system is
well-posed.

\subsection{Steady geostrophic modes} 
On the $f$-plane (planar domain with constant $f$), geostrophic
balanced states satisfying $f\MM{u}^\perp+c^2\nabla\eta=0$ are steady
since $\nabla\cdot\MM{u}=0$. The remaining solutions of the linear
rotating shallow-water equations are fast inertia-gravity waves. In
the quasi-geostrophic limit (slow, large scale motion), when nonlinear
terms and spatially varying $f$ are introduced, these steady states
become slowly-evolving balanced states that characterise large-scale
weather systems. It is crucial that a discretisation gives rise to
steady geostrophic states on the $f$-plane, otherwise when nonlinear
terms and spherical geometry are introduced, balanced states will emit
noisy inertia-gravity waves that will pollute the numerical solution
over timescales that are much shorter than that required for a weather
forecast. To show that mixed finite element methods have steady
geostrophic modes, we follow the approach of \citet{ThRiSkKl2009},
namely we aim to show that vanishing divergence implies steady
vorticity, then checking that vanishing divergence and steady
vorticity implies steady velocity.

To obtain a geostrophic balanced state corresponding to a given
streamfunction $\psi^h$, we initialise $\MM{u}^h$ and $\eta^h$ as
follows:
\begin{enumerate}
\item Set $\MM{u}^h=\nabla^\perp\psi^h$.
\item Set $\eta^h$ from the geostrophic balance relation
\begin{equation}
\label{psi projection}
c^2 \int_{\Omega}\alpha^h\eta^h\diff{V} =
f\int_\Omega \alpha^h\psi^h\diff{V},
\quad
\forall \alpha^h\in V.
\end{equation}
\end{enumerate}
Substitution in equation \eqref{fe u eqn} then gives
\begin{eqnarray*}
\dd{}{t}\int_\Omega \MM{w}^h\cdot\MM{u}^h\diff{V} & = & 
-f\int_\Omega \MM{w}^h\cdot \nabla \psi^h\diff{V} - c^2 
\int_\Omega \nabla\cdot \MM{w}^h \eta^h\diff{V}, \\
& = & f\int_\Omega\nabla\cdot\MM{w}^h \psi^h\diff{V} - c^2 
\int_\Omega \nabla\cdot \MM{w}^h \eta^h\diff{V}, \\
& = & 0,
\end{eqnarray*}
having noted that $\nabla\cdot\MM{w}^h\in V$ and so we may choose
$\alpha^h=\nabla\cdot\MM{w}^h$ in equation \eqref{psi projection}. To
show that $\dot{\eta}^h=0$, first note that
$\MM{u}^h=\nabla^\perp\psi^h$ and hence
$\nabla\cdot\MM{u}^h=0$. Equation \eqref{fe eta eqn} thus becomes
\[
\int_{\Omega} \alpha^h\dot{\eta}^h\diff{V}=0, \quad \forall \alpha^h\in V,
\]
and hence $\dot{\eta}^h=0$. This means that the geostrophic balanced
state is steady.

 \subsection{Numerical dispersion relations}
\label{dispersion relations}
 In this section we consider the numerical wave propagation properties
 of this family of finite element discretisations, on the $f$-plane and
 on the $\beta$-plane in the quasi-geostrophic limit.

 Dispersion relations are computed by assuming time-harmonic solutions
 proportional to $e^{-i\omega t}$ (a valid assumption if the equations
 are invariant under time translations) and studying the resulting
 eigenvalue problem. For the continuous equations on the periodic
 plane, the equations are also invariant under spatial translations
 and so it may be assumed that the eigensolutions take the form
 $e^{i(\MM{k}\cdot\MM{x}-\omega t)}$ where $\MM{k}$ is restricted so
 that the periodic boundary conditions are satisfied. Substitution in
 the equations of motion leads to an algebraic system relating
 $\MM{k}$ to $\omega$: the dispersion relation. For the linear
 shallow-water equations this system is most easily obtained by using
 the Helmholtz decomposition for $\MM{u}$. Numerical dispersion
 relations for continuous-time spatial discretisations are also
 computed by assuming time-harmonic solutions, leading to a discrete
 eigenvalue problem. If a structured mesh is used on the periodic
 plane with a set of discrete translation symmetries then
 eigensolutions satisfy the property that translating from one cell to
 another by $\Delta\MM{x}$ results in the discrete eigensolution
 changing by a factor of $e^{i(\MM{k}\cdot\Delta\MM{x})}$, where
 $\MM{k}$ is again chosen so that the periodic boundary conditions are
 satisfied. This can again lead to a numerical relationship between
 $\MM{k}$ and $\omega$, obtained for both the $f$-plane, and the
 $\beta$-plane in the quasi-geostrophic limit, for the hexagonal
 C-grid in \citet{Th08}, and for the $P1_{DG}-P2$ finite element pair
 in \citet{CoHa2011}.

 Here, we discuss the properties of the discrete eigenvalue problem
 arising from the finite element spaces from the framework of this
 paper. The discussion makes use of the discrete Helmholtz
 decomposition. In the $f$-plane case, substitution of the discrete
 Helmholtz decomposition into equations (\ref{fe u eqn}-\ref{fe eta
   eqn}) and assuming time-harmonic solutions yields
 \begin{eqnarray}
   -i\omega\int_{\Omega}\nabla\gamma^h\cdot\nabla\psi^h\diff{V} +
   \int_{\Omega}f\nabla\gamma^h\cdot D\phi^h\diff{V} & = & 0, 
   \label{fplane psi eqn}
   \\
   \label{fplane phi eqn}
   -i\omega\int_{\Omega} D\alpha^h\cdot D\phi^h\diff{V} + 
   \int_{\Omega}fD\alpha^h\cdot \left(\nabla\psi^h
     +\left(D\phi^h\right)^\perp\right)\diff{V} -
   c^2\int_{\Omega}\nabla\cdot D\alpha^h\eta^h\diff{V} &=& 0, \\
   -i\omega\int_{\Omega}\alpha^h\eta^h\diff{V} +
   \int_{\Omega}\alpha^h\nabla\cdot D\phi^h\diff{V}&=&0,
 \label{fplane eta eqn}
 \end{eqnarray}
 for all test functions $\alpha^h\in V$, $\gamma^h\in E$. Next we 
 define projections $P^E:V\to E$ and $P^V:E\to V$ by
\begin{eqnarray*}
\int_\Omega \nabla\gamma^h\cdot\nabla \left(P^E\phi^h\right)\diff{V} & = &
\int_\Omega \nabla\gamma^h\cdot D\phi^h\diff{V}, \qquad \forall
\phi^h \in V, \quad \gamma^h \in E, 
\\
\int_\Omega D\alpha^h\cdot D \left(P^V\psi^h\right)\diff{V} & = &
\int_\Omega D\alpha^h\cdot \nabla\psi^h\diff{V}, \qquad \forall
\psi^h \in E, \quad \alpha^h \in V.
\end{eqnarray*}
These projections are uniquely defined since $P^E$ uses the standard
continuous finite element discretisation of the Laplace operator which
is solvable by the Lax-Milgram theorem when $E$ is restricted to mean
zero functions, and $P^V$ uses the mixed finite element discretisation
of the Laplace operator using the spaces $S$ and $V$ which is solvable
by the stability conditions \eqref{inf sup} and \eqref{S1} when $V$ is
also restricted to mean zero functions.

Using these projections, and the fact that the divergence operator
maps from $S$ to $V$, equations (\ref{fplane psi eqn}-\ref{fplane eta
  eqn}) become
 \begin{eqnarray}
   -i\omega\psi^h + fP^E\phi^h & = & 0, \\
\nonumber
   -i\omega\int_{\Omega} D\alpha^h\cdot D\phi^h\diff{V} + 
   f\int_{\Omega}D\alpha^h\cdot D P^V\psi^h\diff{V} & & \\
\qquad\qquad     +   \int_{\Omega}fD\alpha^h\cdot\left(D\phi^h\right)^\perp\diff{V} -
   c^2\int_{\Omega}\nabla\cdot D\alpha^h\eta^h\diff{V} &=& 0, \\
   -i\omega\eta^h + \nabla\cdot D\phi^h&=&0,
   \label{eta2phi}
 \end{eqnarray}
and elimination of $\psi^h$ and use of the definition of $D$ gives
\begin{eqnarray}
\nonumber
0 & = & 
   \omega\left(\left(\omega^2+f^2\right)\int_{\Omega} \alpha^h\eta^h\diff{V} +
   \int_{\Omega}\alpha^h\eta^h\diff{V} -
   c^2\int_{\Omega}\nabla\cdot D\alpha^h\eta^h\diff{V}\right) \\
& & \quad
   + if^2\int_{\Omega}D\alpha^h\cdot D \left(P^VP^E\phi^h-\phi^h\right)\diff{V} 
     -   \omega\int_{\Omega}fD\alpha^h\cdot\left(D\phi^h\right)^\perp\diff{V},
\label{f disp reln}
\end{eqnarray}
where $\phi^h$ is obtained from equation \eqref{eta2phi}. The first
row of equation \eqref{f disp reln} is the discretisation of the
continuous eigenvalue problem for the rotating shallow-water equations
using the mixed finite element spaces $V$ and $S$. In this case the
eigenvalues of this discrete eigenvalue problem converge to the
eigenvalues of the continuous problem at the optimal rate as described
in \citet{BoBrGa1997}. However, there are two extra terms in the
bottom row of equation \eqref{f disp reln}. The second term converges
to zero for smooth $\phi^h$, and use might be made of spectral
perturbation theory to examine what effect this has on the discrete
eigenvalue problem; we have not yet developed a technique to do this.
However, the impact of the first term in the second row is more
immediately clear, since it involves projecting $\phi^h$ from $V$ to
$E$ and back to $V$ again. If $V$ has larger dimension than $E$, which
is the case for the lowest order Raviart-Thomas element on triangles,
for example, then this double projection will have a kernel, and
$(P^VP^E-1)\phi^h$ will not be small. This leads to spurious branches
of inertia-gravity waves, \emph{i.e.} branches of solutions of the
discrete eigenvalue problem that do not converge to solutions of the
continuous eigenvalue problem as $h\to 0$. See \cite{Da2010} for
numerical examples illustrating this spurious modes, in particular
Figures 2,3 and 4. Hence, $\dim(V)\leq \dim(E)$ is a necessary
condition for the absence of spurious divergent inertia-gravity modes.

A similar approach can be taken to studying the $\beta$-plane
solutions in the quasi-geostrophic limit. 
Substitution of the discrete Helmholtz decomposition into equations
(\ref{fe u eqn}-\ref{fe eta eqn}) and assuming time-harmonic solutions
yields
 \begin{eqnarray}
   -i\omega\int_{\Omega}\nabla\gamma^h\cdot\nabla\psi^h\diff{V} +
   \int_{\Omega}\left(f_0+\beta y\right)\nabla\gamma^h\cdot D\phi^h\diff{V} & = & 0, 
   \label{betaplane psi eqn}
   \\
\nonumber
   -i\omega\int_{\Omega} D\alpha^h\cdot\left(D\phi^h+\nabla^\perp\psi^h\right)\diff{V} + & & \\ 
   \label{betaplane phi eqn}
   \quad \int_{\Omega}\left(f+\beta y\right)D\alpha^h\cdot \left(\nabla\psi^h+\left(D\phi^h\right)^\perp\right)\diff{V} -
   c^2\int_{\Omega}\nabla\cdot D\alpha^h\eta^h\diff{V} &=& 0, \\
   -i\omega\int_{\Omega}\alpha^h\eta^h\diff{V} +
   \int_{\Omega}\alpha^h\nabla\cdot D\phi^h\diff{V}&=&0.
 \label{betaplane eta eqn}
 \end{eqnarray}
 In the usual quasi-geostrophic limit, the leading order solution is
\[
\phi^h_g =0, \quad
\int_{\Omega}f_0D\alpha^h\cdot \nabla\psi^h_g\diff{V}
+c^2\int_{\Omega}\nabla\cdot D\alpha^h\eta^h_g\diff{V} = 0,
\]
where $\phi^h_g$, $\psi^h_g$ and $\eta^h_g$ are the leading order
terms in the low Rossby number expansion of $\phi^h$, $\psi^h$ and
$\eta^h$ respectively. This is the same as the geostrophic steady
state formula for the $f$-plane, and we have
\[
f_0P^V\psi_g^h = c^2\eta_g^h.
\]
The next order in the expansion of
the equations (we do not make use of the next order in the $\phi^h$
equation) is
 \begin{eqnarray}
   \label{ag psi}
   -i\omega\int_{\Omega}\nabla\gamma^h\cdot\nabla\psi^h_g\diff{V} +
   \int_{\Omega}f_0\nabla\gamma^h\cdot D\phi^h_{ag}\diff{V} +
   \int_{\Omega}\beta y\nabla\gamma^h\cdot \nabla^\perp\psi_g^h\diff{V} & = & 0, 
   \\
   \label{ag eta}
   -i\omega\int_{\Omega}\alpha^h\eta^h_g\diff{V} +
   \int_{\Omega}\alpha^h\nabla\cdot D\phi^h_{ag}\diff{V}&=&0.
 \end{eqnarray}
 Again, the embedding property implies that
 $i\omega\eta^h_g=\nabla\cdot D\phi^h_{ag}$. Since $\gamma^h$ is
 continuous and $D\phi^h_{ag}$ has continuous normal components, we
 may integrate by parts in the second two terms in equation
 \eqref{ag psi}, to obtain
\begin{eqnarray*}
0  & = &  -i\omega\int_{\Omega}\nabla\gamma^h\cdot\nabla\psi^h_g\diff{V} -
  i\omega\int_{\Omega}\frac{f_0^2}{c^2}\gamma^h\psi^h_g\diff{V} -
  \int_{\Omega}\beta \gamma^h \pp{}{x}\psi_g^h\diff{V} \\
  & & \qquad +   i\omega\int_{\Omega}\frac{f_0^2}{c^2}\gamma^h\left(1-P^EP^V\right)\psi^h_g\diff{V}.
\end{eqnarray*}
The first line is the continuous finite element approximation to the
Rossby wave eigenvalue problem using the finite element space $E$,
which has convergent eigenvalues. The second line is a perturbation
involving $\left(1-P^EP^V\right)\psi^h_g$ which will not always be
small if $P^EP^V$ has a non-trivial kernel. This will be the case if
$\dim(V)<\dim(E)$, as occurs in the lowest order Brezzi-Douglas-Marini
(BDM1) element on triangles \citep{BrDoMa1985} which has $P1$ as the
streamfunction space, and hence $2\dim(V)=\dim(E)+2-\dim(H)$. If
$P^EP^V$ has a non-trivial kernel, this will lead to spurious Rossby
wave branches of the numerical dispersion relation. We conclude that
$\dim(V)=\dim(E)$ is a necessary condition for avoiding both spurious
divergent modes and spurious irrotational modes. Note that this is not
a sufficient condition since it is still possible for $P^EP^V$ or
$P^VP^E$ to have non-trivial kernel even in this case. This condition
motivates the selection of examples of mixed finite element spaces
given in the next section.

\section{Examples}
\label{examples}

In this section we provide two examples of mixed finite element spaces
that are suitable for constructing pseudo-uniform grids on the sphere,
and that have the additional property that there are exactly twice as
many velocity degrees of freedom as pressure degrees of freedom, which
prevents the presence of spurious mode branches. The first example is
the modified Raviart-Thomas element on quadrilaterals, and the second
example is the Brezzi-Douglas-Fortin-Marini element on triangles.

\subsection{Modified Raviart-Thomas element on quadrilaterals}
There have been several efforts at developing numerical weather
prediction models based on a cubed sphere grid (see \cite{PuLi2007},
for example) in which a grid on the surface of a cube is projected to
a sphere. The drawback in using such is grid is that to obtain a
C-grid finite difference method with stationary geostrophic states,
the scheme of \cite{ThRiSkKl2009} must be used, which requires the
grid to be orthogonal in the sense that lines joining adjacent
pressure nodes must cross cell boundaries at right-angles. On the
cubed sphere, this condition does not produce a pseudo-uniform grid
since elements become clustered near the poles as the resolution is
increased. Mixed finite elements provide extra freedom to design
numerical schemes since the orthogonality condition is not a
requirement; it is replaced by the conditions on finite element spaces
specified in Section \ref{framework}. 

The lowest-order Raviart-Thomas finite element space is the mixed
finite element analogue of the C-grid since the pressure space is
piecewise constant functions, and the velocity fields are constrained
to be have constant, continuous normal components on element
edge. This means that one normal component of velocity must be stored
on each element edge, just like the C-grid. The velocity fields are
constructed on a square $1\times 1$ reference element $\hat{K}$ with
coordinates $(\xi_1,\xi_2)$, on which the $\xi_1$-component of
velocity $\hat{\MM{u}}$ is obtained by linear interpolation between
constant values on the $\xi_1=0$ and $\xi_1=1$ edges, and the
$\xi_2$-component is obtained by linear interpolation between constant
values on the $\xi_2=0$ and $\xi_2=1$ edges. In these coordinates, the
divergence is constant. In any physical element $K$ in the mesh, we
define a coordinate mapping $\MM{g}:\MM{\xi}\mapsto\MM{x}$, and the
velocity in $K$ is obtained \emph{via} the Piola transformation
\[
\MM{u}(\MM{x}) = \frac{1}{\det\left(\pp{\MM{g}}{\MM{\xi}}\right)}
\pp{\MM{g}}{\MM{\xi}}
\cdot\hat{u}(\MM{\xi}),
\]
which preserves flux integrals 
\[
\int_\gamma\hat{\MM{u}}\cdot\MM{n}\diff{S}(\MM{\xi}) = 
\int_{\MM{g}(\gamma)}\MM{u}\cdot\MM{n}\diff{S}(\MM{x}),
\]
guaranteeing continuity of normal fluxes. The divergence satisfies
\[
\nabla\cdot\MM{u} = \frac{1}{\det\left(\pp{\MM{g}}{\MM{x}}\right)}
\hat{\nabla}\cdot\hat{\MM{u}},
\]
where $\hat{\nabla}$ is the divergence in the local coordinates
$\MM{\xi}$. If the coordinate transformation is affine (elements are
parallelograms), the determinant of the Jacobian is constant, and so
the divergence of the velocity is constant in each element. However,
for general quadrilateral elements (required for the cubed sphere),
the coordinate transformation is bilinear, with linear determinant of
the Jacobian. The solution, proposed by \cite{BoGa2009}, is to modify
the basis functions by adding a divergent correction with vanishing
normal components on the boundary that makes the divergence constant.
The corresponding streamfunction space $E$ is the usual continuous
bilinear space on quadrilaterals, often denoted $Q1$, and it can
easily be shown that the $\nabla^\perp$ operator maps from $E$ into
$S$ in this case. In fact, the Boffi-Gastaldi correction adds a purely
divergent component to the velocity field and so the $\nabla^\perp$
embedding property is not affected.

The RT0-Q0 finite element space has one pressure degree of freedom per
quadrilateral element, and one velocity degree of freedom per
edge. Since (for periodic boundary conditions or the sphere) each edge
is shared by two elements, this means that there are exactly twice as
many velocity degrees of freedom as pressure degrees of freedom. This
modified Raviart-Thomas finite element space satisfies all the
conditions that we require in this paper and hence has potential for
use on pseudo-uniform grids for numerical weather prediction.

\subsection{Brezzi-Douglas-Fortin-Marini element on triangles}
\label{BDFM1}
There is an analogous Raviart-Thomas finite element space on triangles
which satisfies the required embedding properties. However, these
spaces satisfy $2\dim(V)>\dim(S)$ in general. For example, the lowest
order finite element space RT0-P0 has one pressure degree of freedom
per element, and one velocity degree of freedom per edge, meaning that
$3\dim(V)=2\dim(S)$. The BDM1 element on triangles has one pressure
degree of freedom per element and two velocity degrees of freedom per
edge, meaning that $3\dim(V)=\dim(S)$, so $2\dim(V)<\dim(S)$. However,
the little-used lowest order Brezzi-Douglas-Fortin-Marini (BDFM1)
element together with $P1_{DG}$ on triangles satisfies
$2\dim(V)=\dim(S)$. The BDFM family of elements for quadrilaterals was
introduced in \cite{Br+1987}, and an analogous family for triangles
was described in \cite{BrFo1991}. On triangles it is infrequently used
since the BDM and RT families have less degrees of freedom for the
same order of convergence (after suitable post-processing). However,
these extra degrees of freedom are useful to us here since they mean
that $\dim(V)=\dim(E)$.

Here we describe the BDFM1 element on triangles as an augmentation of
the BDM1 element on triangles, which we recall first. Given a triangle
$K$, we define $P_k(K)$ to be the space of $k$-th order polynomials on
$K$. We define the following spaces on $K$:
\begin{eqnarray*}
\mbox{velocity space} &\quad& S(K) = \{P_1(K)\}^2 \\
\mbox{pressure space} &\quad& V(K) = P_0(K).
\end{eqnarray*}
For a triangulation $T$ of the domain $\Omega$, we define the $BDM1$
velocity space
\[
S = \{\MM{v}\in H(\mathrm{div},\Omega):\MM{v}|_K\in S(K), \quad K\in T\},
\]
where $H(\mathrm{div},\Omega)$ is the space of vector fields with
square integrable divergence, which requires that $\MM{v}$ has
continuous normal component across triangle edges. The pressure space
is
\[
V = \{\eta:\eta|_K\in V(K)\},
\]
with no continuity requirements across edges. 

\begin{figure}
\centerline{\includegraphics[width=8cm]{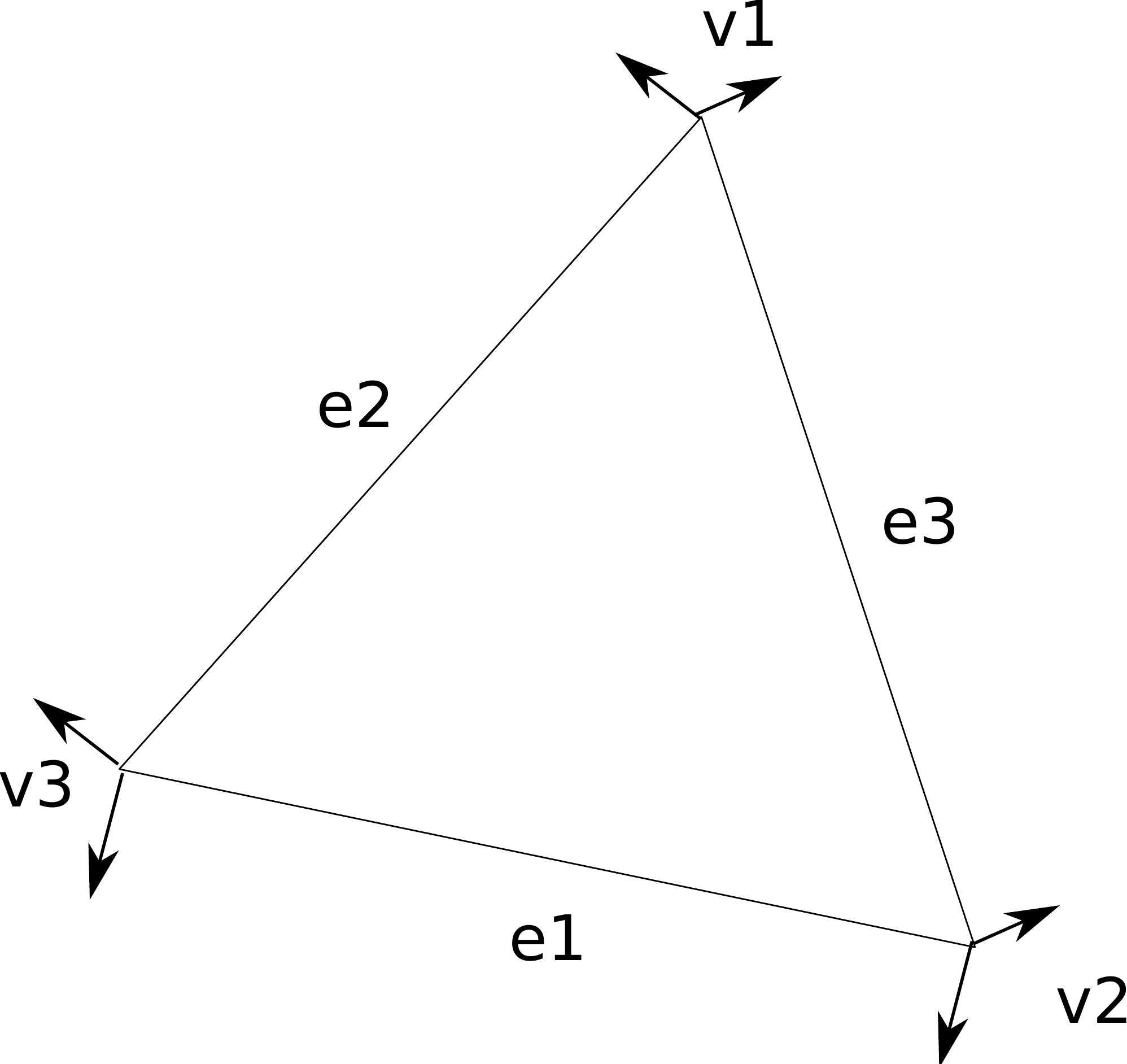}\includegraphics[width=8cm]{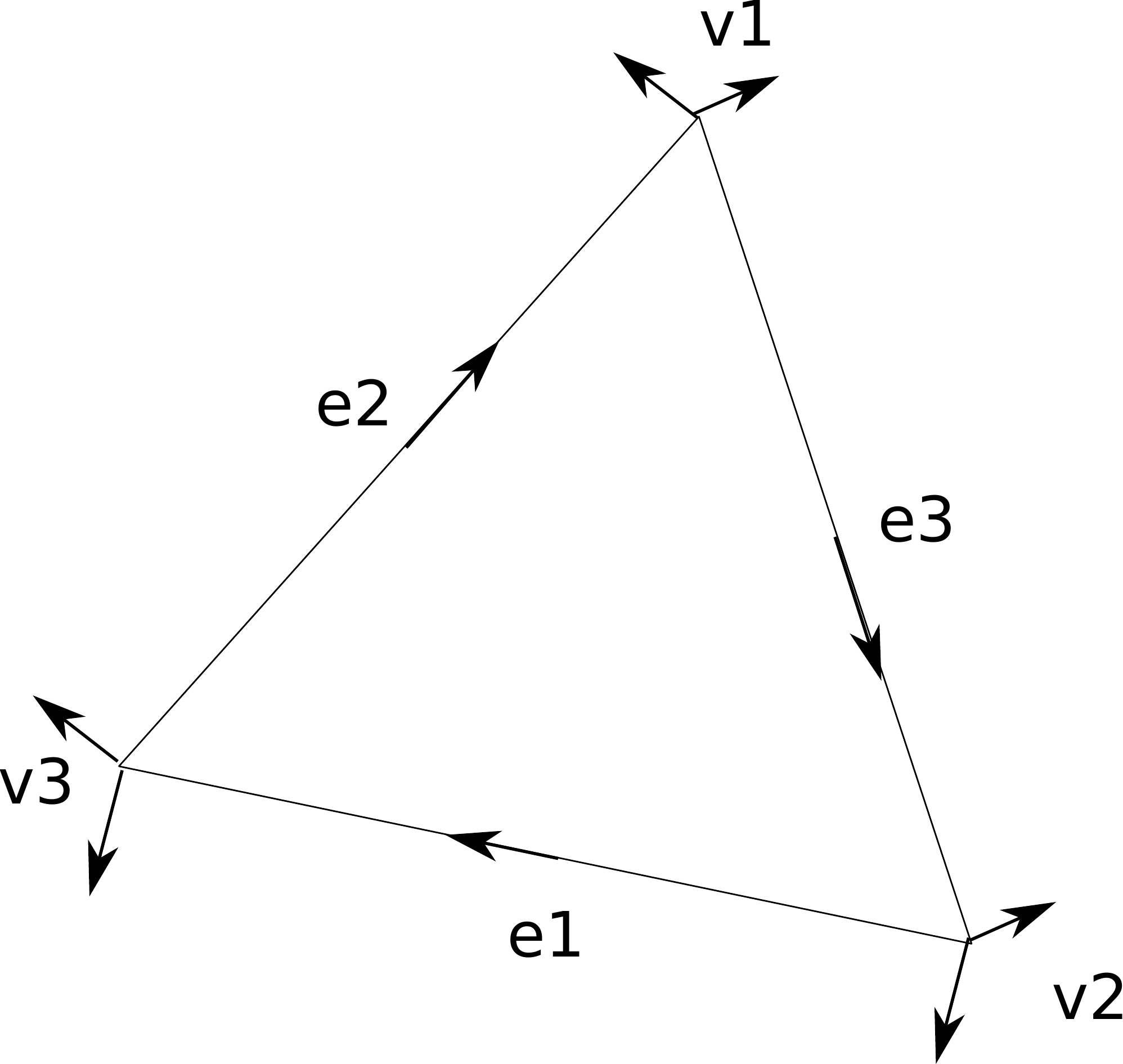}}
\caption{\label{bbm1}Diagram showing degrees of freedom in (left) BBM1
vector element, (right) augmented BBM1 vector element.}
\end{figure}

A convenient set of local nodal basis functions for $S$ is defined by
choosing two node points on each triangle edge, each node located at
one of the vertices belonging to that edge: a total of six node
points.  For example, in the triangle shown in Figure \ref{bbm1}, on
edge $e1$ there are two node points, one at vertex $v3$ and one at
vertex $v2$. The basis function associated with edge $e1$ and vertex
$v3$ is
\[
\MM{\phi}_{1,3} = \MM{t}_2\lambda_3,
\]
where $\MM{t}_2$ is the unit tangent vector to edge $e2$ and where
$\{\lambda_i\}_{i=1}^3$ are the barycentric coordinates associated
with vertices $e1$, $e2$ and $e3$ respectively. It can easily be
checked that $\MM{\phi}_{1,3}$ has normal component equal to 1 at the
node point located at vertex $v3$ on edge $e1$, and normal component
equal to zero at each of the other node points. The other six basis
functions are constructed in a similar manner.

To increase the number of degrees of freedom in each triangle $K$ in
the triangulation $T$, we define the local BDFM1 space $\hat{S}(K)$ by
\[
\hat{S}(K) = \{\MM{v}\in P_2(K)^2:\MM{v}\cdot\MM{n}=0 \mbox{ on } 
\partial K\}.
\]
Since all of the vectors in $\hat{S}(K)$ vanish on the boundary of
$K$, they do not alter the values of the normal components at the
boundary, and so there are no additional continuity constraints.  The
dimension of $P_2(K)^2$ is 12, and there are 9 independent degrees of
freedom which do not vanish on the boundary, which means that
$\dim(\hat{S}(K))=3$.

A convenient set of local nodal basis functions for $\hat{S}$ is
defined by locating nodes that store the tangential component of
velocity at the centre of each edge. The tangential component of
velocity is permitted to be discontinuous and so a different value of
the tangential component will be stored on each side of the edge. The
basis function associated with the node at the centre of edge $e1$ is
\[
\MM{\phi}_1 = 4\MM{t}_1\lambda_2\lambda_3.
\]
It can easily be checked that $\MM{\phi}_1$ has vanishing normal
component on all edges, tangential component equal to 1 at the centre
of edge $e1$ and vanishing tangential component on the other two
edges. The other two basis functions are constructed in a similar manner.

The augmented velocity space $S$ on the triangulation $T$ is defined
as 
\[
S = \{\MM{v}\in
H(\mathrm{div},\Omega):\MM{v}=\MM{v}'+\hat{\MM{v}},
\MM{v}'|_K\in S(K),\quad\hat{\MM{v}}|_K\in\hat{S}(K),\quad K\in T\}.
\]
The pressure space $V$ is defined as
\[
V=\{\eta\in P_1(K)\}
\]
with no continuity requirements. For this mixed element pair the
velocity space $S$ has 6 DOF per element, and the pressure space $V$
has 3 DOF per element, hence there are twice as many velocity DOF as
pressure DOF, just as for the C-grid finite difference method on quadrilaterals.

For our augmented velocity space, it is easy to define the projection
operator $\Pi^S$. The projection is computed element by element and
guarantees the continuity of $\MM{u}\cdot\MM{n}$ across element edges.
The projection on an element $K$ is defined from the following
conditions:
\begin{eqnarray}
\label{BDM1 proj}
\int_{e(i)} \gamma^h(\Pi^S\MM{u}-\MM{u})\cdot\MM{n}\diff{S} & = & 0
\quad \forall \gamma^h\in P^1(e(i)),
\forall\,\mathrm{edges}\,e(i)\in\partial K, i=1,2,3, \\
\label{Shat divergent}
\int_{K}\nabla\gamma^h\cdot(\Pi^S\MM{u}-\MM{u})\diff{V} & =& 0
\quad \forall\gamma^h \in P^1(K), \\
\label{Shat rotational}
\int_{K}\nabla^\perp B\cdot(\Pi^S\MM{u}-\MM{u})\diff{V} & =& 0,
\end{eqnarray}
where $B$ is the cubic ``bubble'' function (as used in the MINI
element \citep{ArBrFo1984}). In a triangle $K$, the cubic bubble
function $B_K$ is the unique cubic polynomial which takes the value 1
at the barycentre and 0 on all three edges. The streamfunction space
$E$ is
\[
E = \{\psi\in H^1(\Omega):\psi|_K=\psi'|_K+\alpha B_K, \psi'|_K\in
P_2(K),\alpha\in \mathbb{R}\}.
\]
Equation \eqref{BDM1 proj} comprises the BDM1 projection operator,
fixing six degrees of freedom. The components of the extra degrees of
freedom $\hat{S}(K)$ are not affected since they all satisfy
$\MM{u}\cdot\MM{n}=0$ on $\partial K$. The vector field $\nabla^\perp
B$ lies inside $\hat{S}(K)$ since it is quadratic (being the skew
gradient of a cubic function, $B$) and has vanishing normal component
on $\partial K$ (since $B$ is zero on $\partial K$). If we construct
an orthogonal (relative to the $L_2$ inner product) decomposition of
$\hat{S}(K)$ into $\nabla^\perp B \oplus \tilde{S}(K)$ then we see
that equation \eqref{Shat rotational} only involves the $\nabla^\perp
B$ component, and equation \eqref{Shat divergent} only involves the
remaining two $\tilde{S}(K)$ components, as
\[
\int_K\nabla\gamma^h\cdot\nabla^\perp B\diff{V} =
-\int_K\underbrace{\nabla^\perp\cdot\nabla\gamma^h}_{=0}B\diff{V}
+\int_{\partial K}\nabla\gamma^h\cdot\MM{n}\underbrace{B}_{=0}\diff{S},
\]
because $B$ vanishes on $\partial K$. The space
$\{\MM{v}=\nabla\gamma^h,\gamma^h\in P1(K)\}$ is spanned by constant
vector fields, and hence equation \eqref{Shat divergent} fixes the two
degrees of freedom in $\tilde{S}(K)$. Bounds on $\Pi^S$ can be
obtained by following the steps of \citet{BrDoMa1985}, since it simply
involves $L_2$ projection onto various moments.

We define the streamfunction space $E$ as the usual Lagrange
continuous quadratic space augmented by cubic bubble functions.  For
any function $\psi\in E$, the curl $\nabla^\perp$ maps into $S$:
$\nabla^\perp\psi\in S$. Furthermore, we may define a projection
operator $\Pi^E:H^1(\Omega)\to H(\mathrm{div})$ by
\begin{eqnarray*}
\Pi^E\psi(v_i) & = & \psi(v_i)\,\forall\,\mathrm{vertices}\,v_i,\,i=1,2,3, \\
\int_{e_i}\Pi^E\psi\diff{S} & = & \int_{e_i}\psi\diff{S},
\,\forall\,\mathrm{edges}\,e_i\,i=1,2,3, \\
\int_K\Pi^E\psi\diff{V} & = & \int_K\psi\diff{V},
\end{eqnarray*}
for each element $K$. To show that the projections commute with
$\nabla^\perp$, \emph{i.e.}
$\Pi^S\nabla^\perp\psi=\nabla^\perp\Pi^E\psi$, we check each of the
conditions (\ref{BDM1 proj}-\ref{Shat rotational}). Condition
\eqref{BDM1 proj} becomes
\begin{eqnarray}
\nonumber
\int_{e(i)}\gamma^h\nabla^\perp\psi\cdot\MM{n}\diff{S} & = & 
\int_{e(i)}\gamma^h\nabla\psi\cdot\diff{\MM{x}}, \\
\nonumber
& = & -\int_{e(i)}\psi\nabla\gamma^h\diff{\MM{x}}
+ [\gamma^h\psi]_{v^-_{e(i)}}^{v^+_{e(i)}}, \\
\nonumber
& = & -\int_{e(i)}\Pi^E\psi\nabla\gamma^h\diff{\MM{x}}
+ [\gamma^h\Pi^E\psi]_{v^-_{e(i)}}^{v^+_{e(i)}}, \\
& = & \int_{e(i)}\gamma^h\nabla^\perp\Pi^E\psi\cdot\MM{n}\diff{S}, \quad
\forall \gamma^h\in P^1(e(i)), \quad i=1,2,3,
\label{BBM1 vort proj}
\end{eqnarray}
where $v^{\pm}_{e(i)}$ are the two vertices at either end of edge
$e(i)$, and having noted that $\nabla\gamma^h$ is constant for
$\gamma^h\in P^1(e(i))$. Condition \eqref{Shat divergent} becomes
\begin{eqnarray*}
\int_{K}\nabla\gamma^h\cdot\Pi^S\nabla^\perp\psi\diff{V} & =&
\int_{K}\nabla\gamma^h\cdot\nabla^\perp\psi\diff{V}, \\
& =&
-\int_{K}\gamma^h\underbrace{\nabla\cdot\nabla^\perp\psi}_{=0}\diff{V}
+ \int_{\partial K}\gamma^h\nabla^\perp\psi\cdot\MM{n}\diff{S}, \\
& =&
\int_{K}\nabla\gamma^h\cdot\nabla^\perp\Pi^E\psi\diff{V}, 
\quad \forall\gamma^h \in P^1(K), \\
\end{eqnarray*}
where we have used equation \eqref{BBM1 vort proj}. Finally, condition
\eqref{Shat rotational} becomes
\begin{eqnarray*}
\int_{K}\nabla^\perp B\cdot\Pi^S\nabla^\perp\psi\diff{V} 
& =& \int_{K}\nabla^\perp B\cdot\nabla^\perp\psi\diff{V}  \\
& =& -\int_{K}\nabla^2 B\psi\diff{V}  + \int_{\partial K}
\underbrace{\nabla^\perp B\cdot\MM{n}}_{=0}\psi\diff{S}, \\
& =& -\int_{K}\nabla^2 B\Pi^E\psi\diff{V}, \\
& =& \int_{K}\nabla^\perp B\cdot\nabla^\perp \Pi^E\psi\diff{V}, 
\end{eqnarray*}
since $\nabla^2 B$ is constant in $K$ and $B$ is zero on $\partial K$.


Counting global degrees of freedom,
\[
\dim(E)=N_{\edge}+N_{\verts}+N_{\face} = 2N_{\edge}+C, \quad
\dim(S)=3N_{\edge}, \quad \dim(V)=2N_{\edge}+3N_{\face},
\]
where $C$ is the Euler characteristic of the domain $\Omega$ which is
equal to 0 for the doubly-periodic domain and equal to 2 on the
sphere. On the sphere there are two extra constraints: namely that the
divergence and the vorticity both integrate to zero, and so in both
cases $\dim(E)+\dim(V)=\dim(S)$.  Finally, we note that each triangle
has three edges which are each shared with one other triangle, and
hence $2N_{\edge}=3N_{\face}$.

\section{Numerical results}
\label{numerics}
In this section we illustrate the properties of the BDFM1 finite
element space applied to the linear rotating shallow-water equations.
The equations were integrated numerically using the implicit midpoint
rule, and the resulting discrete system was solved by using
hybridisation which is a standard technique for solving elliptic
problems (see \cite{BrFo1991} for a detailed description) in which the
continuity constraints on the velocity space are dropped, and are
instead enforced in the equation by Lagrange multipliers. It becomes
possible to eliminate both the velocity and free surface variables
from the matrix equation, leaving a symmetric positive definite system
to solve for the Lagrange multipliers. The velocity and free surface
variables can then be reconstructed element-by-element. One of the
benefits of this approach is that it can be applied when the Coriolis
term is present, resulting in a fully implicit treatment of this term.
In our numerical tests this system was solved using a direct solver.
In the case of BDFM1-\pdg, there are three Lagrange multipliers per
element.

In the test cases with variable Coriolis parameter $f$, a continuous
piecewise quadratic representation of $f$ was used.

\subsection{Steady states for the $f$-plane}
\label{steady states}
We verified that the geostrophic states are exactly steady on the
$f$-plane for the BDFM1 finite element space by randomly generating
streamfunction fields $\psi$ from the streamfunction space $S$ on the
same mesh as used for the $P1_{DG}-P2$ finite element pair steady
state tests in \citet{CoHaPa2009}, with streamfunction equal to zero
on the boundary. This mesh is a planar unstructured mesh in the $x-y$
plane in a $1\times 1$ square region. The velocity was initialised by
setting $\MM{u}=\MM{k}\times\nabla\psi$ where $\MM{k}$ is the unit
normal to the domain \emph{i.e.} $\MM{k}=(0,0,1)$, and $\eta$ was obtained
by solving the discrete elliptic system
\begin{eqnarray}
\int_\Omega \MM{w}^h\cdot\MM{v}^h\diff{V} + 
\int_\Omega c^2\nabla\cdot\MM{w}^h \eta^h\diff{V} & = & 0 \\
\int_\Omega \alpha^h\nabla\cdot\MM{v}\diff{V} & = & 
\int_\Omega D\alpha^h\cdot f\left(\MM{u}^h\right)^{\perp}\diff{V},
\end{eqnarray}
with $c^2=f=1$. We then integrated the equations forward for arbitrary
lengths of time and observed that the layer thickness $h$ and velocity
$\MM{u}$ remained constant up to machine precision. We also conducted
the same experiment on an icosehedral mesh of the unit sphere with
$c^2=f=1$ (following the ``$f$-sphere'' experiment of
\citet{ThRiSkKl2009}) and obtained the same result.

\subsection{Kelvin waves in a circular basin}
Coastal Kelvin waves provide a challenging test since they propagate
at the gravity wave speed along the coast but are geostrophically
balanced in the direction normal to the coast. We used the Kelvin wave
initial condition for a circular basin with unit dimensionless radius
as proposed in \citet{HaKrStPi2007}, with $\Ro=0.1$ and $\Fr=1$. We
integrated the equations until 10 dimensionless time units with a time
step size $\Delta t=0.01$.

\begin{figure}
\centerline{\includegraphics[width=12cm]{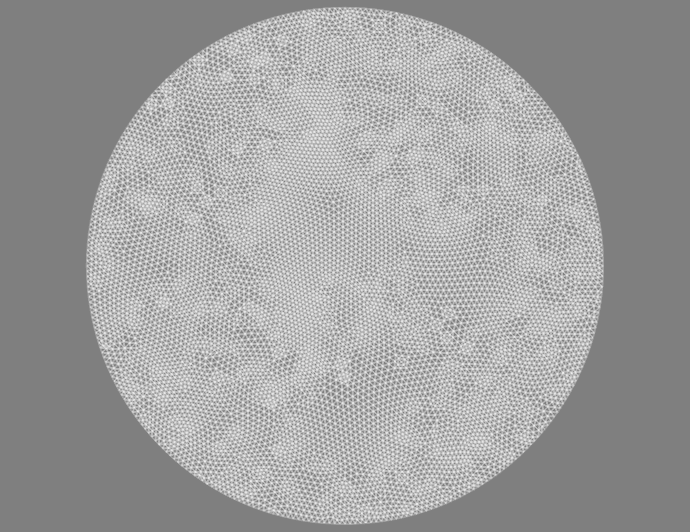}}
\caption{\label{kelvin mesh}Mesh used for the Kelvin wave tests.}
\end{figure}

\begin{figure}
\centerline{\includegraphics[width=6cm]{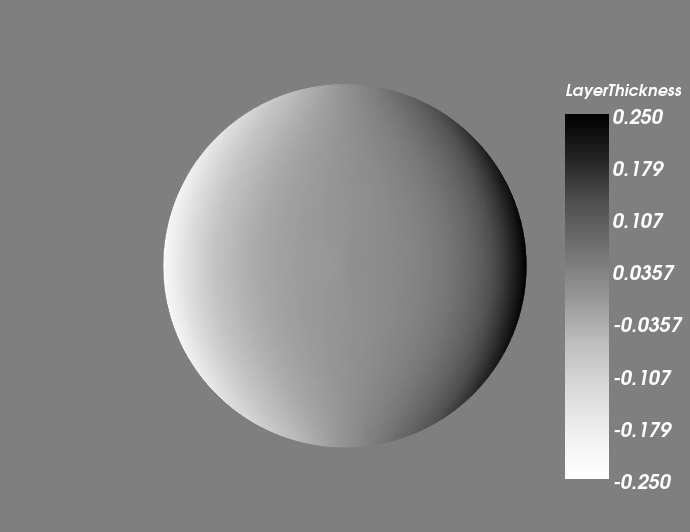}
\includegraphics[width=6cm]{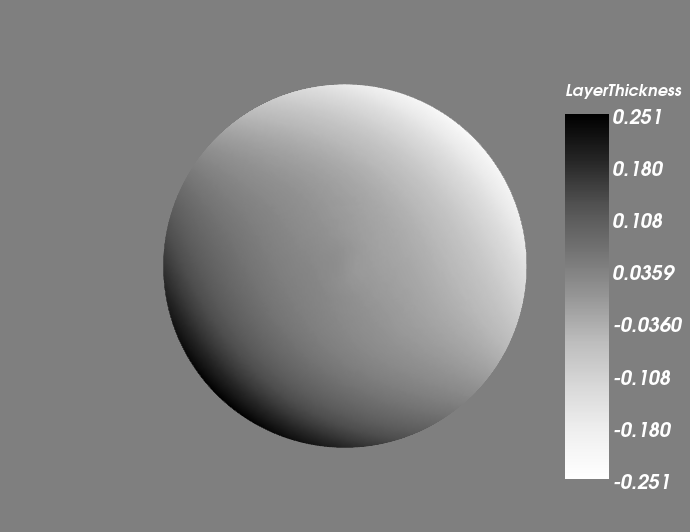}}
\centerline{\includegraphics[width=6cm]{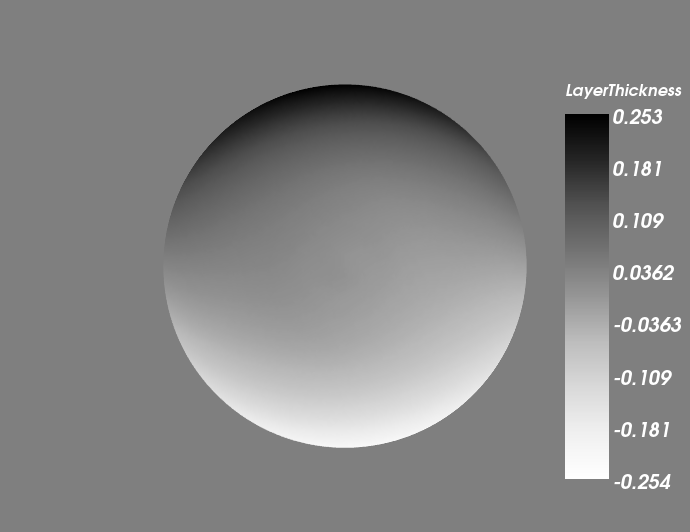}}
\caption{\label{kelvin snapshots}Snapshots of the free surface
  elevation for the circular Kelvin wave testcase obtained at times
  $t=0,2500000,5000000$. The numerical scheme maintains the
  geostrophic balance in the normal direction, as indicated by the
  lack of radiated inertia-gravity waves.}
\end{figure}

The mesh used for the Kelvin wave calculation is shown in Figure
\ref{kelvin mesh}. Some snapshots of the numerical solution are shown
in Figure \ref{kelvin snapshots}. There are no spurious gravity waves
observed, which means that the BDFM1 discretisation is maintaining
geostrophic balance in the normal direction as well as the Kelvin wave
structure.

\subsection{Rossby waves}
\begin{figure}
\centerline{\includegraphics[width=12cm]{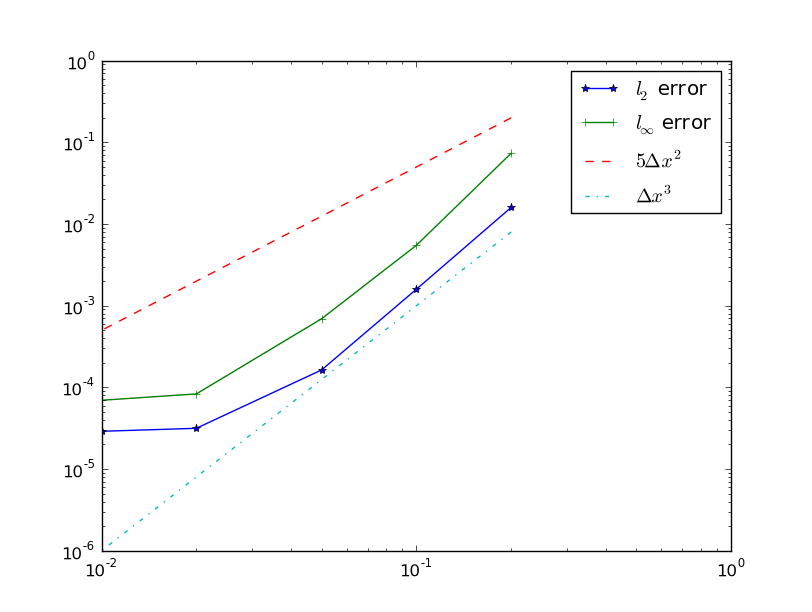}}
\caption{\label{rossby convergence fig}Plot of errors from the Rossby
  convergence test with Rossby number $\Ro=1e-3$ and timestep size
  $\Delta t=0.007996$. The comparison is made after time
  $\pi/(1+8\pi^2)/2$ after which time the wave has travelled halfway
  around the domain. For large $\Delta x$ we observe third-order
  convergence in both $l_2$ and $l_\infty$ norms; for smaller $\Delta
  x$ the error is dominated by either the timestepping error or the
  $\mathcal{O}(\Ro^2)$ truncation error in the small Rossby number
  expansion.}
\end{figure}

To verify the convergence of the method we compared against 
the Rossby wave solution with streamfunction
\[
\psi(x,y,t) = \sin(2\pi x)\sin\left(2\pi\left(y+\gamma t\right)\right),
\qquad \gamma = \frac{2\pi}{1+8\pi^2},
\]
in a square domain with nondimensional length 1, with nondimensional
wave propagation speed $c=Ro^2$, and non-dimensional Coriolis parameter
\[
f = \frac{1 + \Ro y}{\Ro},
\]
and periodic boundary conditions in the $x$-direction. This is an
exact solution of the Rossby wave equation, but is only an asymptotic
limit solution of the linearised rotating shallow-water equations as
$\Ro\to 0$, with $\mathcal{O}(\Ro^2)$ error. This means that for
sufficiently small grid width and time step size we expect the
$\mathcal{O}(\Ro^2)$ error to dominate. The numerical solution was
initialised from this streamfunction following the balanced
initialisation approach described in Section \ref{steady states}. A
plot of the error is shown in figure \ref{rossby convergence fig}. We
observe $\mathcal{O}(\Delta x^3)$ convergence until the error
saturates because of the finite Rossby number. We attribute this third
order convergence to the fact that in Section \ref{dispersion
  relations} the discrete Rossby wave equation was shown to be equal
to usual continuous finite element discretisation of the Rossby wave
equation using the space $E$, plus a perturbation. Since $E$ contains
all of the continuous piecewise quadratic functions, we would expect
third-order convergence provided that the perturbation converges to
zero sufficiently fast (although we do not currently have any
estimates for the convergence of the perturbation).

\begin{figure}
\centerline{\includegraphics[width=5cm]{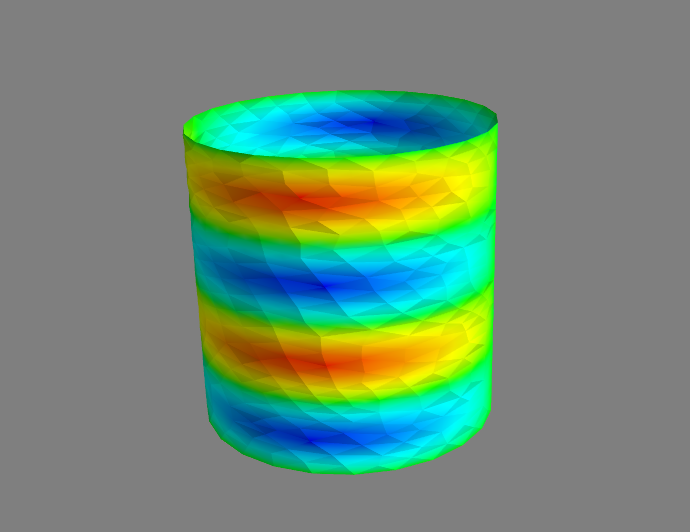}
\includegraphics[width=5cm]{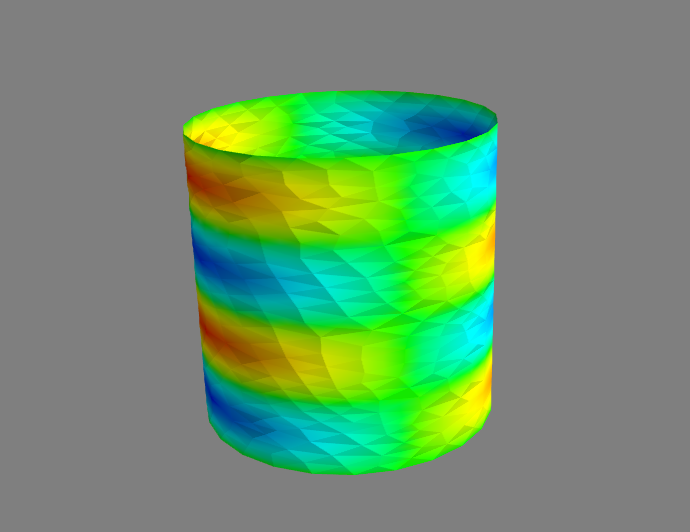}
\includegraphics[width=5cm]{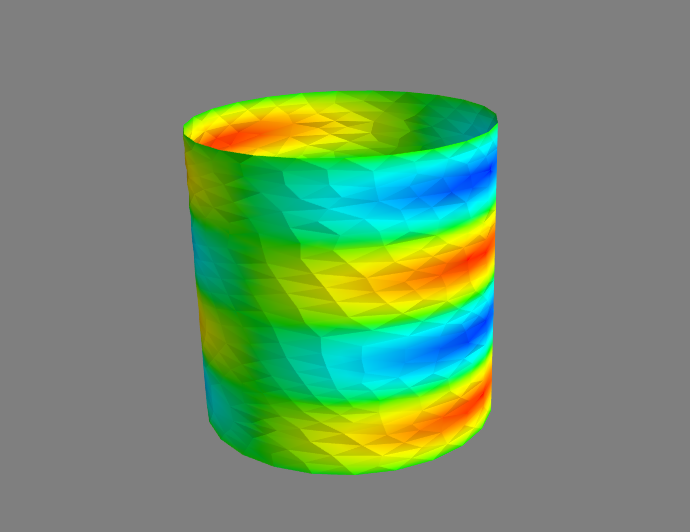}}
\centerline{
\includegraphics[width=5cm]{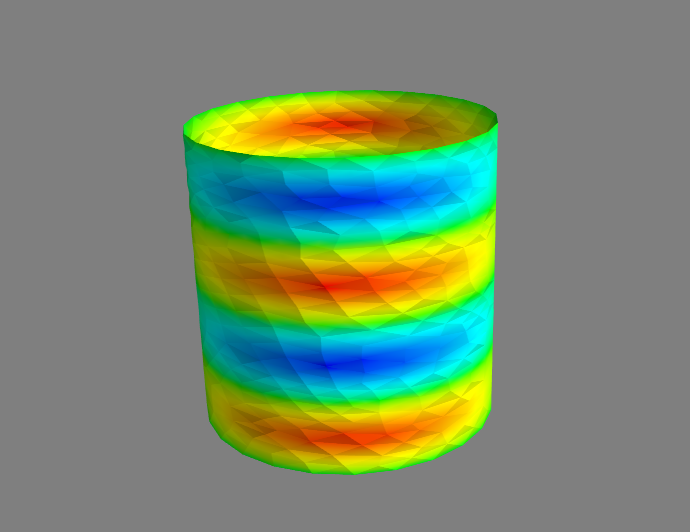}
\includegraphics[width=5cm]{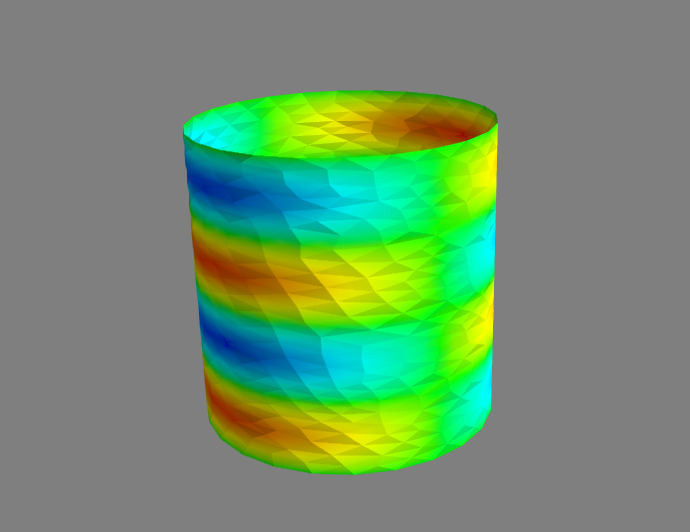}
\includegraphics[width=1.35cm]{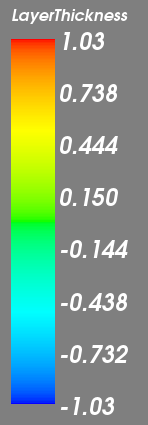}
}
\caption{\label{rossbytube}Rossby waves on the ``$\beta$-tube''
  initialised from a streamfunction on a cylinder with a coarse
  unstructured triangle mesh. Colour plots of the free surface
  elevation are plotted at non-dimensional times
  $t=0.79957,19.9892,39.9784,59.9686,79.9568$ from left to right. No
  unbalanced motions are visible from the plot.}
\end{figure}

To demonstrate the performance of the numerical scheme on arbitrary
manifolds we constructed an unstructured mesh of a cylinder with unit
dimensionless radius and dimensionless height equal to 2. The Coriolis
parameter was set to $f = (1+\Ro z)/\Ro$ and other parameters were
kept the same as the planar Rossby wave tests. We call this
configuration the ``$\beta$''-tube since it corresponds to a
$\beta$-plane that has been wrapped into a cylinder. Some plots of the
numerical integration of this test case are provided in Figure
\ref{rossbytube}; no unbalanced motions are visible from the plots.

\subsection{Solid rotation on the sphere}
To investigate the grid imprinting caused by the finite element
scheme, we integrated the linear rotating shallow-water equations on
the sphere with initial condition obtained from the streamfunction
$\psi=-u_0\cos\theta$, where $\theta$ is the latitude, $u_0=2\pi
R/(12\mbox{ days})$, and $R=6.37122\times 10^6$ is the radius of the
sphere. The rotation rate $|\Omega|$ was $1/(1\mbox{ day})$, and
$g=9.8$. This solution is a steady state solution of the linear
equations with varying $f$ because of the cylindrical symmetry; in
general we do not expect numerical discretisations which break this
symmetry to preserve the steady state.

In our experiment, we used a level 4 icosahedral mesh (each
icosahedron edge being subdivided into 8) of the sphere. The velocity
and free surface elevation were initialised according to the procedure
described in Section \ref{steady states}. To measure the deviation
from a steady state, the free surface elevation after 10 days of
simulation with a timestep of 3600s was subtracted from the initial
condition.  Remarkably, the errors were almost indistinguishable from
round-off error. It turns out that this is because of the mapping used
between functions on the sphere, and functions on the icosahedral mesh
with flat triangular elements used for the numerical integration. The
finite element streamfunction $\psi^h$ was initialised according
to $\psi^h=\psi\circ\phi$, where $\phi$ is the mapping given as
follows:
\[
\phi(x,y,z) = \left(
\left(\frac{R^2-z^2}{x^2+y^2}\right)^{1/2}x,
\left(\frac{R^2-z^2}{x^2+y^2}\right)^{1/2}y,z
\right).
\]
This mapping preserves the value of $z$, and rescales $x$ and $y$ onto
the sphere. Hence, we obtain $\psi^h=z$, which can be represented
exactly in the streamfunction space $E$. The same mapping is also
applied to the finite element representation $f^h$ of the
Coriolis parameter $f$, and we obtain $f^h=2|\Omega|z$ which can
also be represented exactly. Following the balanced initialisation
procedure, the finite element free surface elevation field
$\eta^h$ is obtained by projecting the mapping
$\eta\circ\phi^{-1}$ into the pressure space $V$, where $\eta$ is the
continuous balanced free surface elevation. Substitution into the 
velocity equation gives
\begin{eqnarray*}
  \dd{}{t}\int_{\Omega}\MM{w}^h\cdot\MM{u}^h\diff{V} & = & 
  -\int_{\Omega} f^h \MM{w}^h\cdot\left(\MM{u}^h\right)^\perp
  \diff{V} + c^2\int_{\Omega}\nabla\cdot\MM{w}^h \eta^h\diff{V}, \\
  \mbox{[definition of $f^h$, $\psi^h$ and $\eta^h$]} & = & 
  \int_{\Omega} f \MM{w}^h\cdot\nabla\psi\diff{V}
  + c^2\int_{\Omega}\nabla\cdot\MM{w}^h \eta\diff{V}, \\
  \mbox{[integration by parts]} & = & 
  \int_{\Omega} \MM{w}^h\cdot\nabla\left(\underbrace{f\psi-c^2\eta}_{=0}\right)\diff{V} = 0,
\end{eqnarray*}
where the second step follows since $\nabla\cdot\MM{w}^h\in V$ and so
we can use the fact that $\eta^h$ is a finite element projection of
$\eta$ in $V$, and where in the last step integration by parts was
possible since $\eta$ is continuous and $\MM{w}^h$ has continuous
normal components.

\begin{figure}
\centerline{\includegraphics[width=8cm]{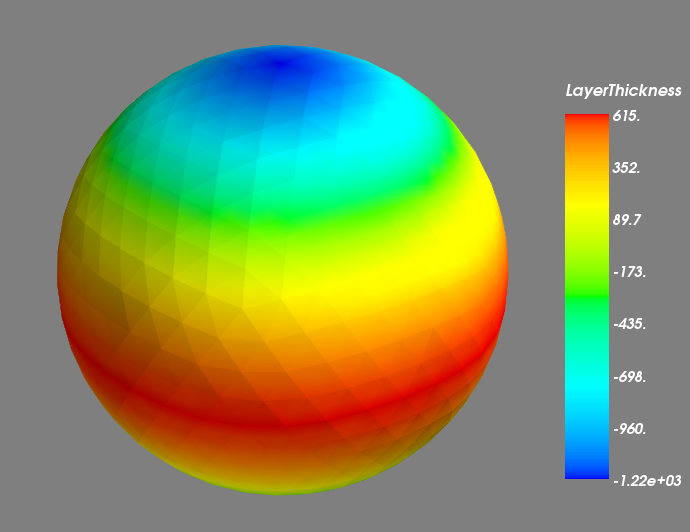}\\
  \includegraphics[width=8cm]{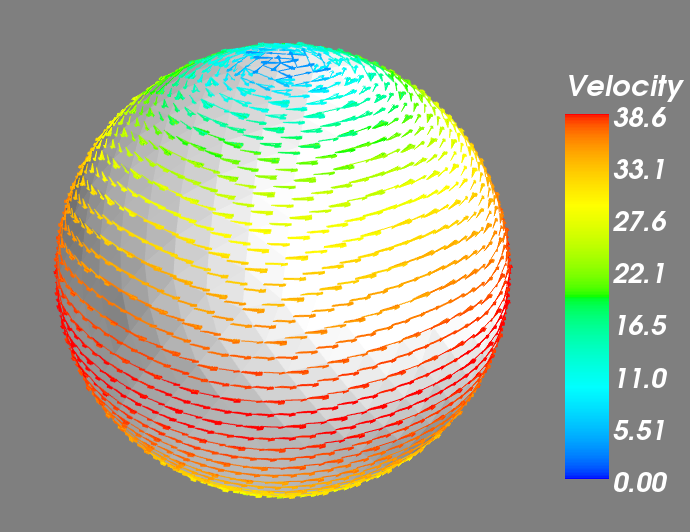}}
\centerline{\includegraphics[width=6cm]{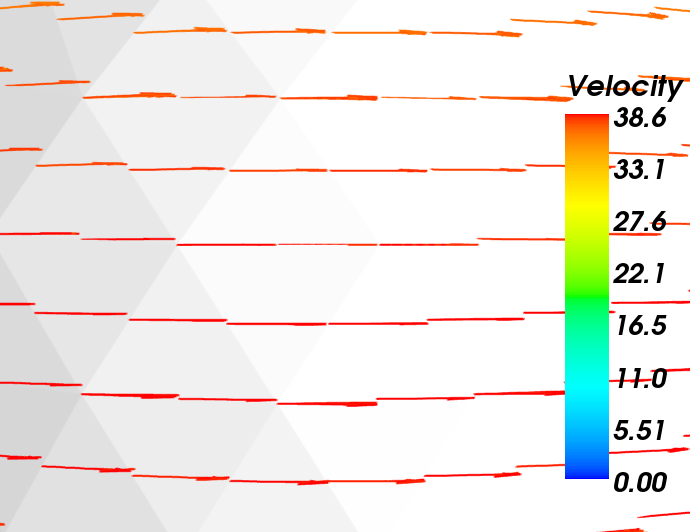}}
\caption{\label{exact steady}Plots showing the exact steady
  numerical solution obtained using the balanced initialisation
  procedure. Top Left: The free surface elevation field. Top Right:
  The velocity field, plotted by evaluating the finite element field
  at vertices and edge midpoints of each triangle. Since only the
  normal components are continuous, there are multiple values of these
  vectors corresponding to the different elements that share those
  vertices/midpoints. Bottom: Close-up of the velocity vectors near
  the equator.}
\end{figure}



\section{Summary and outlook}
\label{summary}

In this paper we described some properties of applying finite element
spaces satisfying the div and curl embedding properties, applied to
the rotating linear shallow water equations, in order to illustrate
their possible suitability for numerical weather prediction on
quasi-uniform grids. In this context, these methods can be thought of
as more flexible extensions of the mimetic C-grid finite difference
method that is currently used in many dynamical cores. This extra
flexibility means that non-orthogonal grids and grids with rapid
changes of mesh resolution can be used, and the ratio of pressure and
velocity degrees of freedom can be adjusted to avoid spurious mode
branches. We showed that spurious inertia-gravity mode branches will
exist if $\dim(E)<\dim(V)$ and spurious Rossby mode branches will
exist if $\dim(V)>\dim(E)$. The discrete Helmholtz decomposition
implies that $\dim(E) = \dim(S)-\dim(V)+2-\dim(H)$ where $H$ is the
space of harmonic velocity fields on the chosen domain. This motivates
the search for finite element spaces with $\dim(S)=2\dim(V)$ that can
be used on pseudo-uniform grids on the sphere. In Section
\ref{examples} we gave two low-order examples: the modified RT0-Q0
element pair for the non-orthogonal cubed sphere, and the BDFM1-\pdg
element pair for triangles, the latter of which was illustrated with
some numerical examples in Section \ref{numerics}.

 In future work, we shall aim to benchmark the augmented mixed element
 pair in the context of numerical weather prediction and ocean
 modelling. One of the benefits of this pair is that discontinuous
 Galerkin methods can be used for the nonlinear continuity equation for
 the density. These methods are locally conservative, have minimal
 dispersion and diffusion errors, and can be made TVB in combination
 with appropriate slope limiters as described in
 \cite{CoSh2001}. Furthermore, as described in \citep{WhLeDe2008}, if
 one wishes to have tracer transport that is both conservative and
 consistent, it is necessary use the pressure space for tracers
 too. This means that tracer transport can (must) also use the
 discontinuous Galerkin method. 
 
 Finally, we note that although the BDFM(k)-Pk$_{DG}$ finite element
 spaces do not have a 2:1 ratio of velocity DOFs to pressure DOFs,
 there does exists a family of higher-order versions of the BDFM1
 element pair with a 2:1 ratio, obtained by appropriately augmenting
 the $\mathrm{BDM}(k)$ spaces (with $k>0$ odd) with higher-order
 components that vanish on element boundaries. This does not work out
 so neatly for $k>1$ since it is also necessary to augment the $P(k)$
 space for pressure, to obtain stable element pairs with twice as many
 velocity DOF as pressure DOF per triangle. In future work, we shall
 investigate these higher-order element pairs, as well as extensions
 to tetrahedra in three-dimensions that can be used in unstructured
 mesh ocean modelling.

\section{Acknowledgements}
The authors acknowledge funding from NERC grants NE/I016007/1,
NE/I02013X/1, and NE/I000747/1. The shallow water solver was produced
using the Imperial College Ocean Model finite element library,
unstructured meshes were generated using GMSH, and the icosahedral
mesh generator was provided by John Thuburn. Plots were obtained using
the Python Matplotlib library and Mayavi2.

\bibliographystyle{elsarticle-harv}
\bibliography{balancedelement}

\end{document}